\documentclass[10pt,a4paper,twoside,notitlepage]{article}

\usepackage{amsmath}
\usepackage{amssymb}
\usepackage{mathrsfs}

\newtheorem{thm}{Theorem}[section]

\newtheorem{prop}{Proposition}[section]
\newtheorem{defn}{Definition}[section]
\newtheorem{rem}{Remark}[section]

\newcommand{\R}{\mathbb{R}}
\newcommand{\N}{\mathbb{N}}
\newcommand{\C}{\mathbb{C}}
\newcommand{\ve}{\varepsilon}
\newcommand{\n}{\noindent}
\newcommand{\vatt}{\mathbb{E}}

\newcommand{\Div}{\mathrm{div}}
\newcommand{\Tr}{\mathrm{tr}}
\newcommand{\D}{\mathrm{D}}
\newcommand{\Id}{\mathrm{I}}
\newcommand{\Img}{\mathrm{Im}}
\newcommand{\re}{\mathrm{Re}}
\newcommand{\Mat}{\mathrm{Mat}}
\newcommand{\Sym}{\mathrm{Sym_d}}
\newcommand{\Symp}{\mathrm{Sym_d^+}}

\newcommand{\prob}{\mathscr{P}}

\begin{document}

\title{Linear-Quadratic $N$-person and 
Mean-Field Games: Infinite Horizon with Discounted Cost and Singular Limits}
 
\author{Fabio S. Priuli
\thanks{Istituto per le Applicazioni del Calcolo ``M. Picone'', C.N.R., Rome, Italy. e-mail: \texttt{f.priuli@iac.cnr.it}
  }}
             
\date{\today}

\maketitle

\begin{abstract}
We consider stochastic differential games with $N$ nearly identical players, linear-Gaussian dynamics, and infinite horizon discounted quadratic cost. Admissible controls are feedbacks for which the system is ergodic.
We first study the existence of affine Nash equilibria by means of an associated system of $N$ Hamilton-Jacobi-Bellman and $N$ Kolmogorov-Fokker-Planck partial differential equations, proving that for small discount factors quadratic-Gaussian solutions exist and are unique. Then, we prove the convergence of such solutions to the unique quadratic-Gaussian solution of the pair of Mean Field equations. We also discuss some singular limits, such as vanishing noise, cheap control and vanishing discount.\\

\n{\bf Key words:} linear-quadratic differential games, Nash equilibria, mean field games.\\

\n{\bf AMS subject classifications:} 49N70, 49N10, 49L25, 91A23

\end{abstract}

\section{Introduction}\label{sec:intro}

\n In this paper we consider an $N$--person differential game driven by a stochastic system of differential equations
\begin{equation}\label{eq:sde}
dX^i_t=(AX^i_t-\alpha^i_t)dt+\sigma dW^i_t\,, \qquad X_0^i=x^i\in\R^d\,,\qquad i=1,\ldots,N\,,
\end{equation}
where $A,\sigma$ are given $d\times d$ matrices, with $\det(\sigma)\neq 0$, $(W^1_t,\ldots,W^N_t)$ are $N$ independent Brownian motions and each $\alpha^i_t\colon [0,+\infty[\to\R^d$ is a bounded process adapted to $W^i_t$ which represents the control of the $i$--th player. Each player wants to minimize on the infinite time horizon a discounted quadratic cost functional given by
\begin{equation}\label{eq:cost_disc}
J^i(X,\alpha^1,\ldots,\alpha^N):=\vatt\left[\int_0^{+\infty}e^{-\ell t}\left(\,{(\alpha_t^i)^TR\alpha_t^i\over 2}\,+f^i(X^i_t\,;m^1,\ldots,m^N)
\right)\,dt\right]
\end{equation}
where $X=(x^1,\ldots,x^N)\in\R^{N d}$ is the initial position of the dynamics, $\vatt$ denotes the expected value, $\ell$ is a positive discount factor, $R$ is a positive definite symmetric $d\times d$ matrix, and we set
\begin{equation}\label{eq:average}
f^i(x;m^1,\ldots,m^N):=
\int_{\R^{d(N-1)}} F^i(\xi^1,\ldots,\xi^{i-1},x,\xi^{i+1},\ldots\xi^N)\prod_{j\neq i} dm^j(\xi^j)\,,
\end{equation}
with
\begin{align}\label{eq:quad_cost}
F^i(X^1,\ldots, X^N)&:= (X-\overline{X_i})^T Q^i(X-\overline{X_i})
=\sum_{j,k=1}^N (X^j-\overline{X_i}^j)^T Q^i_{jk}(X^k-\overline{X_i}^k)\,,
\end{align}
for suitable  $Nd\times Nd$ symmetric matrices $Q^i$ and suitable reference positions $\overline{X_i}\in\R^{N d}$. The notation $Q^i_{jk}$ ($j,k\in\{1,\ldots,N\}$) is used for the $d\times d$ block matrices of $Q^i$.
In~\eqref{eq:cost_disc} and~\eqref{eq:average}, we denoted with $m^1,\ldots, m^N$ the invariant measures associated to the processes $X^1_t,\ldots, X^N_t$. In other words, we are assuming that the cost $J^i$ depends directly on the state of the $i$--th player only, while the other players only influence the cost through their asymptotic distribution in the environment, since $f^i$ represents an average of the quadratic cost $F^i$ w.r.t. the invariant measures of other players. 
The standing assumptions on the game~\eqref{eq:sde}--\eqref{eq:cost_disc} are summed up in the following conditions.
\begin{description}
\item{\bf (H1)} The matrix $\sigma$ in~\eqref{eq:sde} is invertible, the matrix $R$ is symmetric and positive definite and the matrices $Q^i$ in~\eqref{eq:quad_cost} are symmetric.
\item{\bf (H2)} There exist matrices $Q,B$, $C_1,\ldots,C_N$, $D_1,\ldots,D_N$ and vectors $\Delta, H$ such that block matrices and reference states in~\eqref{eq:quad_cost} satisfy for all $i$
$$
Q^i_{ii}= Q\in \Symp\,,
\qquad
\overline{X_i}^i=H\,,
$$
$$
Q^i_{ij}=\,{B\over 2}\,,
\qquad
Q^i_{jj}=C_i\,,
\qquad
\overline{X_i}^j=\Delta\,,
\qquad
\qquad
\forall~j\neq i\,,
$$
$$
Q^i_{jk}=D_i\,,
\qquad
\qquad
\forall~j,k\neq i\,,j\neq k\,.
$$
\item{\bf (H3)} The matrix $A$ is symmetric and there exist constants $r,k>0$ such that $R=r\,\Id_d$ and $\nu:=\,{\sigma^T\sigma\over 2}=k\,\Id_d$.
\end{description}
In~\cite{Bardi,BardiPriuli} games satisfying {\bf (H2)} were referred to as games with ``{\em nearly identical players}''. 
Notice that for such games we can rewrite~\eqref{eq:quad_cost} as
\begin{align*}
F^i(X^1,\ldots, X^N)&=(X^i-H)^T\,Q\,(X^i-H)\\
&+\sum_{j\neq i}\left[ (X^i-H)^T\,{B\over 2}\,(X^j-\Delta)+ (X^j-\Delta)^T\,{B\over 2}\,(X^i-H)\right]\\
&+\sum_{j\neq i}^N (X^j-\Delta)^T C_i(X^j-\Delta)+\sum_{j,k\neq i}^N (X^j-\Delta)^T D_i(X^k-\Delta)\,,
\end{align*}
which, in particular, means that each player cannot distinguish among other players and tries to reach his happy state $H$ while pushing all competitors towards a common state $\Delta$.

For games of the form~\eqref{eq:sde}--\eqref{eq:cost_disc}, we study in this paper the existence of Nash equilibria through the solutions of an associated system of $2N$ Hamilton--Jacobi--Bellman (HJB, in the following) and Kolmogorov--Fokker--Planck (KFP) equations
\begin{equation}\label{eq:hjkfp_disc_intro}
\left\{
\begin{array}{l}
-k\Delta v^i+\displaystyle\,{1\over 2r}\,|\nabla v^i|^2-(\nabla v^i)^TAx+\ell v^i=f^i(x;m^1,\ldots,m^N)\\
-k\Delta m^i-\displaystyle\Div\left(m^i \cdot\Big({\nabla v^i\over r}\,-A x\Big)\right)=0\\
\int_{\R^d}m^i(x)\,dx=1\,,\qquad m^i>0
\end{array}
\right.
\qquad\qquad i=1,\ldots, N
\end{equation}
where the unknown $v^i,m^i$ represent respectively the value function for the $i$--th player and its invariant measure (with a slight abuse of notations, we denote with $m^i$ a measure as well as its density), and $\Div$ is the divergence operator.
In view of the Linear-Quadratic structure of the game, we look for solutions of the HJB--KFP system in the class of quadratic value functions and multivariate Gaussian distributions. This produces Nash equilibria for~\eqref{eq:sde}--\eqref{eq:cost_disc} in the form of affine feedbacks. 

Our result for these games is that, for small values of the discount factor $\ell>0$, there exists a unique Quadratic--Gaussian (abbreviated QG later on) solution to~\eqref{eq:hjkfp_disc_intro} and thus a unique affine Nash equilibrium strategy.
Moreover, we rigorously prove that, as the number of players $N$ tends to infinity, QG solutions of~\eqref{eq:hjkfp_disc_intro} converge to solutions of the Mean Field PDE system
$$
\left\{
\begin{array}{l}
-k\Delta v+\displaystyle\,{1\over 2r}\,|\nabla v|^2-\nabla v^TAx+\ell v=\hat V[m](x)\\
-k\Delta m-\displaystyle\Div\left(m \cdot\Big({\nabla v\over r}\,-A x\Big)\right)=0\\
%\displaystyle\textred{\min\left[v(x)-x^T\,{RA\over 2}\,x\right]=0}\,,\qquad
\int_{\R^d}m(x)\,dx=1\,,\qquad m>0
\end{array}
\right.
$$
for a suitable integral operator $\hat V$ mapping probability measures into quadratic polynomials of the variable $x$. This latter result perfectly matches the ones obtained by Lasry \& Lions in their seminal papers~\cite{LL1,LL2,LL3} about differential games on the torus $\mathbb{T}^d$, and the ones on ergodic LQ games in $\R$ and $\R^d$ (see~\cite{Bardi} and~\cite{BardiPriuli}, respectively).

\n Then, we investigate the relation between the games with discounted cost~\eqref{eq:cost_disc} and the ones with long--time--average cost functional studied in~\cite{BardiPriuli}. Namely, using the same notations as above, we consider the ergodic cost
\begin{equation}\label{eq:cost_ergodic}
J^i(X,\alpha^1,\ldots,\alpha^N):= \liminf_{T\to\infty}\,{1\over T}\,\vatt\left[
\int_0^T \,{(\alpha_t^i)^TR\,\alpha_t^i\over 2}\,+F^i(X^1,\ldots, X^N)\,dt
\right]\,,
\end{equation}
whose affine Nash equilibria were characterized in~\cite{BardiPriuli} through the study of the corresponding HJB and KFP equations. Here, we prove that the QG solutions giving Nash equilibria for the game~\eqref{eq:sde}--\eqref{eq:cost_disc} converge as $\ell\to 0^+$ to the corresponding QG solutions for the game~\eqref{eq:sde}--\eqref{eq:cost_ergodic}, as in the case of classical differential games. Moreover, we prove that the limit procedures as $\ell\to 0^+$ and as $N\to+\infty$ commute.

\n Finally, we investigate other singular limits procedures, and prove that the deterministic ($\nu\to 0$) and cheap cost ($R\to 0$) limits for the games with ergodic cost~\eqref{eq:sde}--\eqref{eq:cost_ergodic} do commute with the mean field limit ($N\to+\infty$). 

Linear--Quadratic differential games have a large literature, see the books~\cite{BO,Engw} and the references therein. The Lasry--Lions approach to MFG, originally introduced in~\cite{LL1,LL2,LL3}, has found application to several different contexts spanning from numerical methods~\cite{AcCD}, to discrete games~\cite{Gomes}, to financial problems~\cite{GLL}.
Large population limits for multi--agent systems were also studied independently by Huang, Caines and Malhame~\cite{HCM03,HCM06,HCM07ieee}.
They introduced a method named ``Nash certainty equivalence principle'' that produces a feedback from a mean--field equation, and shows that such control gives an approximate Nash equilibrium for the $N$--person game if $N$ is large enough. We cannot review here the number of papers inspired by their approach, but let us cite~\cite{BSYY,LZ} for LQ problems,~\cite{NCMH} for recent progress on nonlinear systems,~\cite{KTY} on the rate of convergence as $N\to\infty$, and the references therein.
In particular, we mention that~\cite{HCM07ieee,LZ} deal with discounted infinite horizon games as the ones we are considering here. There are some differences between our results and the ones in the cited papers, though. In~\cite{HCM07ieee,LZ} more general costs $J^i$ are allowed, explicitly depending on other players' states $X^j$, but only the existence of approximate Nash equilibria is established. Here, we trade off the generality of the cost to prove existence of \emph{exact} Nash equilibria for the game, and to prove the relation between $N$--players games and their mean field limit, as $N\to+\infty$. More details will be discussed in section~\ref{sec:compare_caines}.

The paper is organized as follows. In section~\ref{sec:prelim} we recall some preliminary facts for symmetric matrices, algebraic Riccati equations and LQ games~\eqref{eq:sde}--\eqref{eq:cost_ergodic}. Section~\ref{sec:discounted} is devoted to the existence of Nash equilibria for infinite horizon differential games with discounted cost~\eqref{eq:sde}--\eqref{eq:cost_disc}. Section~\ref{sec:singular} contains the results about singular limits as $\nu\to 0$ (deterministic limit), $R\to 0$ (cheap control) and $\ell\to 0^+$ (vanishing discount). Finally, section~\ref{sec:tech_proof} contains the proofs of the results and section~\ref{sec:conclusion} discusses extensions and open problems.

\section{Notations and preliminaries}\label{sec:prelim}

\subsection{Matrices and eigenvalues}

In the following, we will use the notation $\Mat_{m\times n}(\R)$ for the linear space of real $m\times n$ matrices, $\Id_d\in\Mat_{d\times d}(\R)$ for the identical $d\times d$ matrix and $\mathrm{spec}(A)$ for the spectrum of a matrix $A$.
The linear subspace of real symmetric $d\times d$ matrices will be denoted by $\Sym$ and, for $M\in\Sym$, we say that $M$ is positive semidefinite (resp. positive definite) if for all $x\in\R^d$ there holds $x^TMx\geq 0$ (resp. if for all $x\in\R^d\setminus\{0\}$ there holds $x^TMx> 0$). The notation $\Symp$ will be used for the set of real symmetric and positive definite $d\times d$ matrices. 
Recall that for matrices $M\in\Sym$, the expression
\begin{equation}\label{eq:max_spec_norm}
\|M\|:=\max\,\{|\ell|~;~\ell\in\mathrm{spec}(M)\}\,,
\end{equation}
defines a norm. In particular, $\|M\|=\max\mathrm{spec}(M)$ whenever $M$ is positive semidefinite. Also, eigenvalues of a matrix depend continuously on its coefficients (see e.g.~\cite{SerreMat}) so that, for instance, given a sequence of symmetric matrices $A_n\to A$, the sequences of the minimal and maximal eigenvalues of $A_n$ converge, respectively, to $\min\,\mathrm{spec}(A)$ and $\max\,\mathrm{spec}(A)$.
We conclude with a property that will be used in the rest of the paper (cf again~\cite{SerreMat}).

\begin{prop}\label{prop:matrices} Let $H\in\Sym$ and $K\in\Symp$. Then, $HK$ is diagonalizable with real eigenvalues and the number of positive (resp. negative) eigenvalues of $HK$ is equal to the number of positive (resp. negative) eigenvalues of $H$. The same holds for $KH$.
\end{prop}

\subsection{Admissible strategies and Nash equilibria}

\begin{defn}\label{defn:admiss_strategy} A strategy $\alpha^i$ is said to be \emph{admissible} (for the $i$--th player) if it is a bounded process adapted to $W^i_t$ such that the corresponding solution $X^i_t$ to~\eqref{eq:sde} satisfies
\begin{itemize}
\item $\vatt[X^i_t]$ and $\vatt[(X^i_t)(X^i_t)^T]$ are both bounded on $[0,T]$ for every $T$;
\item $X^i_t$ is \emph{ergodic} in the following sense: there exists a probability measure $m^i=m^i(\alpha^i)$ on $\R^d$ such that
$$
\int_{\R^d}|x|\,dm^i(x)<\infty
\qquad\qquad
\int_{\R^d}|x|^2\,dm^i(x)<\infty
$$
and
$$
\lim_{T\to+\infty}\,{1\over T}~\vatt\left[\int_0^Tg(X^i_t)\,dt\right]=\int_{\R^d}g(x)\,dm^i(x)\,,
$$
locally uniformly w.r.t. the initial state $X^i_0$, for all functions $g$ which are polynomials of degree at most $2$.
\end{itemize}
\end{defn}

\n In~\cite{BardiPriuli} it was shown that all affine strategies $\alpha^i(x)= K^ix+c^i$ with $K^i\in\Mat_{d\times d}(\R)$ such that the matrix $A-K^i$ has only eigenvalues with negative real part, and $c^i\in\R^d$, are admissible. Namely, considering $\alpha^i_t:= \alpha^i(X^i_t)$ with $X^i_t$ solution of
\begin{equation}\label{eq:csde}
dX^i_t=[(A-K^i)X^i_t-c^i]dt+\sigma^idW^i_t\,,
\end{equation}
$\alpha^i_t$ is admissible and $X^i_t$ has a unique invariant measure $m^i$ given by a multivariate Gaussian.

\begin{defn} A vector of admissible strategies $\overline{\alpha}=(\overline{\alpha}^1,\ldots,\overline{\alpha}^N)$ is a \emph{Nash equilibrium strategy} for the $N$--person game with dynamics~\eqref{eq:sde} and cost $J^i$ given by either~\eqref{eq:cost_ergodic} or~\eqref{eq:cost_disc}, if for every index $i\in\{1,\ldots,N\}$ and for every admissible strategy $\alpha^i$ for the $i$--th player there holds
$$
J^i(X,\overline{\alpha}^1,\ldots,\overline{\alpha}^N)\leq J^i(X,\overline{\alpha}^1,\ldots,\overline{\alpha}^{i-1},\alpha^i,\overline{\alpha}^{i+1},\ldots\overline{\alpha}^N)\,.
$$
The Nash equilibrium is said to be \emph{symmetric} if all the players adopt the same strategy.
\end{defn}

\subsection{Algebraic Riccati equations}

We recall here some basic facts about algebraic Riccati equations (ARE in the following).

\begin{prop}\label{prop:ARE} Consider the ARE
\begin{equation}\label{eq:ARE}
Y {\cal R} Y
 + Y{\cal A}+{\cal A}^T Y
- {\cal Q}=0
\end{equation}
with ${\cal R}\in\Symp$, ${\cal Q}\in\Sym$ and ${\cal A}$ any $d\times d$ matrix, and introduce the following notations
$$
\Xi_S:= \left[
\begin{array}{c}
I_d\\ S
\end{array}
\right]\in\Mat_{2d\times d}(\R)\,,
\qquad\qquad
{\cal H}:=
\left(
\begin{array}{cc}
{\cal A} & {\cal R}\\
{\cal Q} & -{\cal A}^T
\end{array}
\right)\in\Mat_{2d\times 2d}(\R)\,,
$$
where $S$ is any element of $\Mat_{d\times d}(\R)$, and $\Img\,\Xi_S$ for the $d$--dimensional linear subspace of $\R^{2d}$ spanned by the columns of $\Xi_S$. Then the following facts hold.
\begin{description}
\item{\it (i)} $Y$ is a solution of~\eqref{eq:ARE} if and only if $\Img\,\Xi_Y$ is ${\cal H}$--invariant, i.e. if and only if ${\cal H}\xi\in\Img\,\Xi$ for all $\xi\in\Img\,\Xi$.
\item{\it (ii)} If the matrix ${\cal H}$ has no purely imaginary nonzero eigenvalues, then equation~\eqref{eq:ARE} has solutions $Y$ such that $Y=Y^T$.
\item{\it (iii)} If~\eqref{eq:ARE} has symmetric solutions, then there exists a unique symmetric solution $Y$ with
$$
\big\{\lambda\in\mathrm{spec}({\cal A}+{\cal R}Y)~;~\re\,\lambda\neq 0\big\} = \mathrm{spec}({\cal H})\cap \big\{z\in\C~;~\re\, z>0\big\}\,.
$$
In particular, if ${\cal H}$ has only real nonzero eigenvalues, then there exists a unique symmetric solution $Y$ such that
\begin{equation}\label{eq:spec_ARE}
\mathrm{spec}({\cal A}+{\cal R}Y) = \mathrm{spec}({\cal H})\cap (0,+\infty)\,.
\end{equation}
\end{description}
\end{prop}

\n The proof follows from standard arguments about Riccati equations that can be found in~\cite{Engw,LR}. We give here some explicit references for sake of completeness. Part {\it (i)} is contained in Proposition 7.1.1 of~\cite{LR}. Part {\it (ii)} is a particular case of Theorem 8.1.7 in~\cite{LR}. Finally, part {\it (iii)} is proved in Theorem 8.3.2 of~\cite{LR}.

\subsection{Results for LQ games with ergodic cost}\label{sec:LQG}

\n In view of the study of the singular limits, we review the results obtained in~\cite{BardiPriuli} for LQ differential games with ergodic costs. We start by noticing that, 
for the games~\eqref{eq:sde}--\eqref{eq:cost_ergodic}, all players share the same Hamiltonian given by
\begin{align*}
H(x,p)&:=\min_{\omega}\left\{-\omega^T \,{R\over 2}\,\omega - p^T\big(A\,x-\omega\big)\right\}=-p^TA\,x+\min_{\omega}\left\{-\omega^T \,{R\over 2}\,\omega - p^T\cdot\omega\right\}\,.
\end{align*}
Since the minimum is attained at $\omega=R^{-1} p$, we conclude
$H(x,p)=p^T\,{R^{-1}\over 2}\,p - p^TA\,x$.
Therefore, the system of HJB--KFP equations associated to the game is given by
\begin{equation}\label{eq:hjkfp}
\left\{
\begin{array}{l}
-\Tr(\nu\,\D^2 v^i)+\displaystyle\,{1\over 2}\,(\nabla v^i)^T R^{-1}\nabla v^i-(\nabla v^i)^TAx+\lambda^i=f^i(x;m^1,\ldots,m^N)\\
-\Tr(\nu\,\D^2 m^i)-\displaystyle\Div\Big(m^i  \cdot(R^{-1}\nabla v^i-A x)\Big)=0\qquad\qquad\qquad\qquad\qquad i=1,\ldots, N\\
\int_{\R^d}m^i(x)\,dx=1\,,\qquad m^i>0
\end{array}
\right.
\end{equation}
where the unknown $v^i,m^i$ represent respectively the value function for the $i$--th player and its invariant measure, and $\lambda^i$ is a real number representing the outcome of the game for the $i$--th player. Here $\Tr$ and $\Div$ are respectively the trace of a matrix and the divergence operator. In order to formulate the algebraic conditions which characterize the existence of Quadratic--Gaussian (QG in the rest of the paper) solutions to~\eqref{eq:hjkfp}, we need the following definition.

\begin{defn}\label{def:RSprop} Given matrices ${\bf A}\in\Sym$ and ${\bf N}, {\bf R}, {\bf Q}\in\Symp$, we say that $({\bf A}, {\bf N}, {\bf R}, {\bf Q})$ satisfy the \emph{Riccati--Sylvester property} if every symmetric and positive definite solution $Y$ of the ARE
\begin{equation}\label{eq:riccati}
Y\,{{\bf N} {\bf R}{\bf N}\over 2}\,Y=\,{{\bf A}^T {\bf R}{\bf A}\over 2}\,+{\bf Q}\,,
\end{equation}
is also a solution of the Sylvester equation
\begin{equation}\label{eq:sylvester}
Y {\bf N}{\bf R}-{\bf R}{\bf N} Y= {\bf R}{\bf A}-{\bf A}^T{\bf R}\,.
\end{equation}
\end{defn}

The first result for $N$--players games~\eqref{eq:sde}--\eqref{eq:cost_ergodic} satisfying {\bf (H1)} and {\bf (H2)} was the following (cf Theorem~2 in~\cite{BardiPriuli}): The system of $2N$ HJB--KFP equations~\eqref{eq:hjkfp} admits a unique solution $(v^i,m^i,\lambda^i)$ of the form
\begin{equation}\label{eq:ansatz_nearly_id}
v^i(x) =x^T\,{ \Lambda\over 2}\, x+\rho x\,,
\qquad\qquad
m^i(x)= {\cal N}(\mu,\Sigma^{-1})\,,
\qquad\qquad
\lambda^i\in\R\,,
\end{equation}
for suitable symmetric matrices $\Lambda,\Sigma$, with $\Sigma$ positive definite, and suitable vectors $\mu,\rho$, which are in common for all the players, if and only if $(A,\nu,R,Q)$ satisfy the Riccati--Sylvester property in the sense of Definition~\ref{def:RSprop} and the matrix ${\cal B}:= Q+\,{A^TRA\over 2}\,+(N-1)~{B\over 2}$ is invertible. Moreover, the affine feedbacks $\overline{\alpha}^i=\overline{\alpha}:= R^{-1}\nabla v$, for $i=1,\ldots,N$, provide a symmetric Nash equilibrium strategy for all initial positions $X\in\R^{Nd}$ and $J^i(X,\overline{\alpha})=\lambda^i$ for all $X$ and all $i$.

\n In particular, by going through the proof of this Theorem in~\cite{BardiPriuli}, one sees that the coefficients $\Lambda,\Sigma,\rho,\mu$ are determined by solving the following algebraic relations
\begin{equation}\label{eq:KFP_matrix}
\Sigma\,{\nu R\nu\over 2}\,\Sigma-\,{A^T RA\over 2}\,-Q=0\,,
\qquad 
{\cal B} \mu=P\,,
\qquad 
\Lambda=R\big(\nu\Sigma+A\big)\,,
\qquad
\rho=-R\,\nu\,\Sigma\mu\,.
\end{equation}
with $P:=Q H+(N-1)~{B\over 2}\,\Delta$, and $\lambda^i=F^i(\Sigma,\mu)+\Tr(\nu R\nu\Sigma+\nu R A)-\mu^T\,{\Sigma \nu R\nu \Sigma \over 2}\, \mu$
with
\begin{align}\label{eq:v1}
F^i(\Sigma,\mu)&:= H^TQH-(N-1)\,H^T\,{B\over 2}\,(\mu-\Delta)
-(N-1)(\mu-\Delta)^T\,{B\over 2}\,H+(N-1)\Tr(C_i \Sigma^{-1})\nonumber\\
&~~~ + (N-1)(\mu-\Delta)^T C_i(\mu-\Delta)
+(N-1)(N-2)(\mu-\Delta)^T D_i(\mu-\Delta)\,.
\end{align}

\medskip

In order to study the behavior of QG solutions of~\eqref{eq:hjkfp} as $N\to+\infty$, we assume for simplicity that the control system, the costs of the control and the reference positions are always the same, i.e. that $A,\sigma,R,H$ and $\Delta$ are all independent from the number of players $N$. We also denote with
$$
Q^N\,,
\qquad
B^N\,,
\qquad
C_i^N\,,
\qquad
D_i^N\,,
$$
the primary and secondary costs of displacement, respectively, which are assumed to depend on $N$. Concerning these quantities, we require that they tend to suitable matrices $\hat Q,\hat B,\hat C,\hat D$ with their natural scaling, i.e., that as $N\to+\infty$ there hold
\begin{equation}\label{eq:scale}
Q^N\to \hat Q\,,
\qquad
B^N(N-1)\to \hat B\,,
\qquad
C_i^N(N-1)\to \hat C\,,
\qquad
D_i^N(N-1)^2\to \hat D\,,
\qquad
\forall~i\,.
\end{equation}
If we define an operator on probability measures of $\R^d$ by setting for all measures ${\frak m}\in\prob(\R^d)$
\begin{align*}
\hat V[{\frak m}](X)&:=(X-H)^T \hat Q (X-H)
+\!\!\int_{\R^d}\!\!\left((X-H)^T \,{\hat B\over 2}\, (\xi-\Delta)+(\xi-\Delta)^T\,{\hat B\over 2}\,(X-H)\right)d{\frak m}(\xi)\\
&~~~~~+\int_{\R^d} (\xi-\Delta)^T \hat C (\xi-\Delta)\,d{\frak m}(\xi)
+\left(\int_{\R^d} (\xi-\Delta)\,d{\frak m}(\xi)\right)^T\hat D\left(\int_{\R^d} (\xi-\Delta)\,d{\frak m}(\xi)\right)
\end{align*}
then it is possible to verify that, as $N\to+\infty$, the solutions $v^i_N$, $m^i_N$ and $\lambda^i_N$ of~\eqref{eq:hjkfp} tend to solutions of the system of mean field equations
\begin{equation}\label{eq:mfpde}
\left\{
\begin{array}{l}
-\Tr(\nu\D^2 v)+\displaystyle\,{1\over 2}\,\nabla v^T R^{-1}\nabla v-\nabla v^TAx+\lambda=\hat V[m](x)\\
-\Tr(\nu\D^2 m)-\displaystyle\Div\Big(m \cdot(R^{-1}\nabla v-A x)\Big)=0\\
%\displaystyle\textred{\min\left[v(x)-x^T\,{RA\over 2}\,x\right]=0}\,,\qquad
\int_{\R^d}m(x)\,dx=1\,,\qquad m>0
\end{array}
\right.
\end{equation}
like in~\cite{Bardi,LL1,LL3}. Namely, if we assume that
\begin{equation}\label{eq:hyp_limit}
\nu\in\Symp\,,
\qquad\qquad
R\in\Symp\,,
\qquad\qquad
\hat Q\in\Symp\,.
\end{equation}
the following facts hold (cf Theorem~3 in~\cite{BardiPriuli}). First of all, the system~\eqref{eq:mfpde} admits a unique solution $(v,m,\lambda)$ of the form
\begin{equation}\label{eq:ansatz_limit}
v(x) =x^T\,{ \Lambda\over 2}\, x+\rho x\,,
\qquad\qquad
m(x)= {\cal N}(\mu,\Sigma^{-1})\,,
\qquad\qquad
\lambda\in\R\,,
\end{equation}
for suitable symmetric matrices $\Lambda,\Sigma$, with $\Sigma$ positive definite, and suitable vectors $\mu,\rho$,
if and only if $(A,\nu,R,\hat Q)$ satisfy the Riccati--Sylvester property in the sense of Definition~\ref{def:RSprop} and the matrix ${\cal B}^\infty:= \hat Q+\,{A^TRA\over 2}\,+\,{\hat B\over 2}$ is invertible. If in addition $\hat B\geq 0$, then the solution $(v,m,\lambda)$ of the form~\eqref{eq:ansatz_limit} is the unique solution of~\eqref{eq:mfpde} such that  $v(0)=0$. Finally, assume we are given a sequence of $N$--players differential games of the form~\eqref{eq:sde}--\eqref{eq:cost_ergodic} which satisfy {\bf (H1)} and {\bf (H2)} and admit solutions of the form~\eqref{eq:ansatz_nearly_id} for all $N$. Then, if the the limit system~\eqref{eq:mfpde} admits a unique solution of the form~\eqref{eq:ansatz_limit}, we have that the QG solutions $(v^i_N,m^i_N,\lambda_N^i)$ of the $N$--person game converge as $N\to+\infty$ to the QG solution $(v,m,\lambda)$ of~\eqref{eq:mfpde} in the following sense: for all $i=1,\ldots,N$, $v^i_N\to v$ in $C^1_{loc}(\R^d)$ with second derivative converging uniformly in $\R^d$, $m^i_N\to m$ in $C^k(\R^d)$ for all $k$, and $\lambda_N^i\to \lambda$ in $\R$.
%
%Notice that the condition $\hat B\geq 0$ in Theorem~\ref{thm:LIMplay}{\it(b)} can be proven to be equivalent to 
%$$
%\int_{\R^d} \left(\hat V[{\frak m}]-\hat V[{\frak n}]\right)\!(x)\,\, d({\frak m}-{\frak n})(x)\geq 0\,,\qquad\qquad\forall~{\frak m},{\frak n}\in\prob(\R^d)\,,
%$$
%for the integral operator $\hat V$. Thus, it represents a monotonicity property for the operator $\hat V$, along the lines of the uniqueness results in~\cite{Bardi,LL1,LL3}.
%
\n For later use, we also remark that the coefficients $\Lambda,\Sigma,\rho,\mu$ in~\eqref{eq:ansatz_limit} are determined by solving the following algebraic relations
\begin{equation}\label{eq:KFP_matrix_limit}
\Sigma\,{\nu R\nu\over 2}\,\Sigma-\,{A^T RA\over 2}\,-\hat Q=0\,,
\qquad
{\cal B}^\infty \mu=P^\infty\,,
\qquad 
\Lambda=R\big(\nu\Sigma+A\big)\,,
\qquad \rho=-R\,\nu\,\Sigma\mu\,.
\end{equation}
with $P^\infty:=\hat Q H+\,{\hat B\over 2}\,\Delta$, and $\lambda=\hat F(\Sigma,\mu)+\Tr(\nu R\nu\Sigma+\nu R A)-\mu^T\,{\Sigma \nu R\nu \Sigma \over 2}\, \mu$,
with
\begin{align}\label{eq:v1_limit}
\hat F(\Sigma,\mu)&:= H^T\hat Q H-\left(H^T{\hat B\over 2}\,(\mu-\Delta)+(\mu-\Delta)^T {\hat B\over 2}\,H\right)\nonumber\\
&~~~+\Tr(\hat C \Sigma^{-1}) +(\mu-\Delta)^T(\hat C+\hat D)(\mu-\Delta)\,.
\end{align}

\section{Discounted problems} \label{sec:discounted}

In this section, we extend the analysis of~\cite{BardiPriuli} to the case of infinite horizon $N$--person games 
\begin{equation}\label{eq:sde2}
dX^i_t=(AX^i_t-\alpha^i_t)dt+\sigma dW^i_t\,, \qquad X_0^i=x^i\in\R^d\,,\qquad i=1,\ldots,N\,,
\end{equation}
with discounted costs
\begin{equation}\label{eq:cost_disc2}
J^i(X,\alpha^1,\ldots,\alpha^N):=\vatt\left[\int_0^{+\infty}e^{-\ell t}\left(\,{r\,|\alpha_t^i|^2\over 2}\,+f^i(X^i_t\,;m^1,\ldots,m^N)
\right)\,dt\right]\,,
\end{equation}
which satisfy {\bf (H1)}--{\bf (H3)}. In this case, the associated system of $2N$ HJ--KFP equation takes the form
\begin{equation}\label{eq:hjkfp_disc}
\left\{
\begin{array}{l}
-k\Delta v^i+\displaystyle\,{1\over 2r}\,|\nabla v^i|^2-(\nabla v^i)^TAx+\ell v^i=f^i(x;m^1,\ldots,m^N)\\
-k\Delta m^i-\displaystyle\Div\left(m^i \cdot\Big({\nabla v^i\over r}\,-A x\Big)\right)=0\\
\int_{\R^d}m^i(x)\,dx=1\,,\qquad m^i>0
\end{array}
\right.
\qquad\qquad i=1,\ldots, N\,.
\end{equation}

\begin{rem} Observe that if {\bf (H3)} holds, then $(A,\nu,R,Q)$ satisfy the Riccati--Sylvester property. Indeed, equation~\eqref{eq:sylvester} reduces to $k\,r\,(Y-Y)=r(A-A^T)=0$, which is identically satisfied for all $Y\in\Sym$.
\end{rem}

\begin{thm}\label{thm:disc_Nplay} Assume {\bf (H1)}--{\bf (H3)}. Then, there exists $\bar\ell>0$ such that for $\ell<\bar\ell$ the system of HJB--KFP equations~\eqref{eq:hjkfp_disc} admits a unique solution $(v^i_\ell,m^i_\ell)$ satisfying
$$
v^i_\ell(x)=x^T\,{\Lambda_\ell\over 2}\, x+\rho_\ell x+c_\ell^i\,,\qquad\qquad
m^i_\ell(x)= {\cal N}(\mu_\ell,\Sigma_\ell^{-1})\,,
$$
for suitable symmetric matrices $\Lambda_\ell,\Sigma_\ell$, with $\Sigma_\ell$ positive definite, vectors $\mu_\ell,\rho_\ell$ and numbers $c^1_\ell,\ldots,c^N_\ell\in\R$, if and only if the matrix ${\cal B}_\ell:= Q+r\,{A^2\over 2}\,-\ell\,r\,{A\over 2}\,+(N-1)~{B\over 2}$ is invertible. \\
Moreover, the affine feedbacks
$\overline{\alpha}^i(x)=\,{\nabla v^i_\ell(x)\over r}$, for $x\in\R^d$ and $i=1,\ldots,N$,
provide a symmetric Nash equilibrium strategy for~\eqref{eq:sde2}--\eqref{eq:cost_disc2}, for all initial positions $X\in\R^{Nd}$.
\end{thm}

\n The proof is quite technical and it is deferred to section~\ref{sec:tech_proof}. Here we mention that, similarly to the results in section~\ref{sec:LQG}, the coefficients $\Lambda_\ell,\Sigma_\ell,\mu_\ell,\rho_\ell,c^i_\ell$ are characterized by% the following algebraic relations
\begin{equation}\label{eq:KFP_matrix_discount}
\Lambda_\ell=r\big(k\Sigma_\ell+A\big)\,,\qquad\qquad \rho_\ell=-r\,k\,\Sigma_\ell\mu_\ell\,,
\end{equation}
\begin{align}\label{eq:v0_disc}
\ell c^i_\ell&=F^i(\Sigma_\ell,\mu_\ell)+kr\,\Tr(k\Sigma_\ell+A)-\,{k^2r\over 2}\,(\mu_\ell)^T\,\Sigma_\ell^2\, \mu_\ell
\end{align}
where $F^i$ is the function defined by~\eqref{eq:v1} and $\Sigma_\ell$ and $\mu_\ell$ solve respectively
\begin{equation}\label{eq:discount1}
\Sigma_\ell{\cal R}\Sigma_\ell+{\cal A}_{\ell}^T\Sigma_\ell+\Sigma_\ell{\cal A}_{\ell}-{\cal Q}_\ell=0\,,
\qquad\qquad
{\cal B}_\ell \mu_\ell=QH+(N-1)\,{B\over 2}\,\Delta\,.
\end{equation}
with
${\cal R}:=\,{k^2 r\over 2}\,\Id_d$, ${\cal A}_\ell:=\ell\,{k r\over 4}\,\Id_d$ and ${\cal Q}_\ell:=Q+r\,{A^2\over 2}\,-\ell\,r\,{A\over 2}$.

\begin{rem}\label{rem:disc_are} Observe that the conclusion of Theorem~\ref{thm:disc_Nplay} fails when $\ell$ is not small enough. Indeed, the ARE in~\eqref{eq:discount1} may fail to have solutions in $\Symp$, which in turn would give no Gaussian solution $m^i$ for the KFP equation in~\eqref{eq:hjkfp_disc}. To see this, recall that  Proposition~\ref{prop:ARE}{\it (iii)} ensures the existence of a unique solution $Y_\ell\in\Sym$ such that
${\cal A}_\ell+{\cal R}Y_\ell = \,{k^2 r\over 2}\,\left({\ell\over 2 k}\,\Id_d+Y_\ell\right)$
has eigenvalues which coincide with the ones with positive real part of
\begin{equation}\label{eq:discount_are_matrix}
{\cal H}_\ell:=\left(
\begin{array}{cc}
{\cal A}_\ell & {\cal R}\\
{\cal Q}_\ell & -{\cal A}_\ell^T
\end{array}
\right)\,.
\end{equation}
Since $\lambda\in\R$ is an eigenvalue of ${\cal A}_\ell+{\cal R}Y_\ell$ if and only if $\left({2\over k^2 r}\lambda-\,{\ell\over 2 k}\right)$ is an eigenvalue of $Y_\ell$, it is clear that ${\cal A}_\ell+{\cal R}Y_\ell$ has real eigenvalues and that, fixing $\ell>0$ such that $kr>(4\hat\lambda/\ell)$ for some eigenvalue $\hat\lambda> 0$ of ${\cal A}_\ell+{\cal R}Y_\ell$, the matrix $Y_\ell$ has a negative eigenvalue and does not belong to $\Symp$.
Given that no other positive definite solution can exist, because any such $Z_\ell$ would give positive spectrum to ${\cal A}_\ell+{\cal R}Z_\ell$ and this would violate the uniqueness of $Y_\ell$, this means that the game~\eqref{eq:sde2}--\eqref{eq:cost_disc2} corresponding to this value $\ell$ admit no Quadratic--Gaussian solutions to~\eqref{eq:hjkfp_disc}.
\end{rem}

To study the convergence of Nash equilibria as $N\to+\infty$, assume again that the coefficients $A,\sigma,R,H$ and $\Delta$ are all independent from the number of players $N$. Also, assume that the discount factor $\ell$ does not depend on $N$ and that~\eqref{eq:scale} holds for the cost coefficients $Q^N,B^N,C_i^N,D_i^N$.
By denoting with $(v^i_N,m^i_N)$ the solutions found in Theorem~\ref{thm:disc_Nplay}, we expect that they converge, like for games with ergodic costs~\cite{Bardi,BardiPriuli,LL1,LL3}, to solutions of the system of two mean field equations
\begin{equation}\label{eq:mfpde_disc}
\left\{
\begin{array}{l}
-k\Delta v+\displaystyle\,{1\over 2r}\,|\nabla v|^2-\nabla v^TAx+\ell v=\hat V[m](x)\\
-k\Delta m-\displaystyle\Div\left(m \cdot\Big({\nabla v\over r}\,-A x\Big)\right)=0\\
\int_{\R^d}m(x)\,dx=1\,,\qquad m>0
\end{array}
\right.
\end{equation}
Along the lines of Theorem~3 in~\cite{BardiPriuli} (see also section~\ref{sec:LQG}), our main result for this system is the following, whose proof is given in section~\ref{sec:tech_proof}.

\begin{thm}\label{thm:disc_LIMplay} Assume that $r,k>0$ in~\eqref{eq:mfpde_disc} and that the matrix $\hat Q$ in~\eqref{eq:scale} satisfies $\hat Q\in\Symp$. Then, the following facts hold.
\begin{description}
\item{\it (a)}~{\bf [Solutions to MFPDE]} There exists $\hat\ell>0$ such that for $\ell<\hat\ell$ the system~\eqref{eq:mfpde_disc} admits a unique solution $(v,m)$ satisfying
\begin{equation}\label{eq:disc_ansatz_limit}
v(x)=x^T\,{\Lambda\over 2}\, x+\rho x+c\,,\qquad\qquad
m(x)= {\cal N}(\mu,\Sigma^{-1})\,,
\end{equation}
for suitable symmetric matrices $\Lambda,\Sigma$, with $\Sigma$ positive definite, vectors $\mu,\rho$ and $c\in\R$ if and only if the matrix ${\cal B}_\ell^\infty:= \hat Q+r\,{A^2\over 2}\,-\ell\,r\,{A\over 2}\,+\,{\hat B\over 2}$ is invertible. 

\item{\it (b)}~{\bf [Uniqueness]} If in addition $\hat B\geq 0$,  the solution $(v,m)$ of the form~\eqref{eq:disc_ansatz_limit} is the unique solution of~\eqref{eq:mfpde_disc} such that  $v(0)=c$.

\item{\it (c)}~{\bf [Convergence as $N\to\infty$]} Let assume $\ell<\hat\ell$, where $\hat\ell>0$ is the value found in {\it (a)}. For all $N\in\N$ consider $N$--players differential games of the form~\eqref{eq:sde2}--\eqref{eq:cost_disc2} such that {\bf (H1)}--{\bf (H3)} hold. Assume that~\eqref{eq:scale} is verified as $N\to+\infty$,  and that the Mean-Field system~\eqref{eq:mfpde_disc} admits a unique Quadratic--Gaussian solution. Then, the solutions $(v^i_N,m^i_N)$ found in Theorem~\ref{thm:disc_Nplay} converge to a solution $(v,m)$ of~\eqref{eq:mfpde_disc} as $N\to+\infty$ in the following sense: for all $i=1,\ldots,N$, $v^i_N\to v$ in $C^1_{loc}(\R^d)$ with second derivative converging uniformly in $\R^d$ and $m^i_N\to m$ in $C^k(\R^d)$.\\
Moreover such solution is the unique one given in {\it (a)}, with $(v,m)$ of the form~\eqref{eq:disc_ansatz_limit}.
\end{description}
\end{thm}

\section{Singular limits}\label{sec:singular}

We collect in this section, some results on singular limit processes for the LQG $N$--person games and mean field games. 

We start from the result on the vanishing discount limit, which shows the relation between the solutions found in Theorems~\ref{thm:disc_Nplay} and~\ref{thm:disc_LIMplay}, and their limits as the discount factor $\ell$ tends to $0$. We prove that the limit procedures as $\ell\to 0^+$ and as $N\to+\infty$ commute and that both tend to the solution of the mean field equation for the problem with ergodic cost described in section~\ref{sec:LQG}.

\begin{thm}\label{thm:vd} For $N\in\N$, consider $N$--players games of the form~\eqref{eq:sde2}--\eqref{eq:cost_disc2} such that {\bf (H1)}--{\bf (H3)} hold. Assume that~\eqref{eq:scale} holds as $N\to+\infty$ and that ${\cal B}^\infty= \hat Q+\,{A^TRA\over 2}\,+\,{\hat B\over 2}$ is invertible.

\n Then, the vanishing discount limit as $\ell\to 0^+$ and the mean field limit as $N\to+\infty$ commute. Namely, denoting with $(v^i_{\ell,N}, m^i_{\ell,N})$ the solutions to the $N$--players game with discount factor $\ell>0$, there hold
\begin{equation}\label{eq:commute_disc1}
\lim_{\ell\to0^+}\lim_{N\to+\infty}\Big[v^i_{\ell,N}-v^i_{\ell,N}(0)\Big] = \lim_{N\to+\infty}\lim_{\ell\to0^+}\Big[v^i_{\ell,N}-v^i_{\ell,N}(0)\Big] = v
\end{equation}
in  $C^1_{loc}(\R^d)$ with second derivative converging uniformly in $\R^d$,
\begin{equation}\label{eq:commute_disc2}
\lim_{\ell\to0^+}\lim_{N\to+\infty}m^i_{\ell,N} = \lim_{N\to+\infty}\lim_{\ell\to0^+}m^i_{\ell,N} = m
\qquad\qquad \mbox{in  $C^k(\R^d)$ for all $k$,}
\end{equation}
\begin{equation}\label{eq:commute_disc3}
\lim_{\ell\to0^+}\lim_{N\to+\infty}\ell v^i_{\ell,N} = \lim_{N\to+\infty}\lim_{\ell\to0^+}\ell v^i_{\ell,N} = \lambda
\qquad\qquad \mbox{uniformly in $\R^d$,}
\end{equation}
where 
%in the first case the limit is in $C^1_{loc}(\R^d)$ with second derivative converging uniformly in $\R^d$, in the second case the limit is in $C^k(\R^d)$, and 
$(v,m,\lambda)$ is the QG solution to~\eqref{eq:mfpde}.
\end{thm}

\begin{rem} As a byproduct of the previous proof, we have proved that as $\ell\to0^+$ the solution of HJB--KPF system for $N$--players games with discounted cost~\eqref{eq:sde2}--\eqref{eq:cost_disc2} converge to the solution of the corresponding system for $N$--players games with ergodic cost~\eqref{eq:sde}--\eqref{eq:cost_ergodic}. The same holds for solutions of the Mean Field systems of PDE.
\end{rem}

Next we consider the deterministic limit as $k\to 0^+$ (and hence as the noise matrix $\sigma\to 0$) and we prove that such limit and the limit to the Mean Field PDE as $N\to+\infty$ do commute.

\begin{thm}\label{thm:vv} For $N\in\N$, consider $N$--players games of the form~\eqref{eq:sde}--\eqref{eq:cost_ergodic} such that {\bf (H1)}--{\bf (H3)} hold. Assume that~\eqref{eq:scale} holds as $N\to+\infty$ and that ${\cal B}^\infty= \hat Q+\,{A^TRA\over 2}\,+\,{\hat B\over 2}$ is invertible.

\n Then, the deterministic limit as $k\to 0^+$ and the mean field limit as $N\to+\infty$ commute.
Namely, denoting with $(v^i_{k,N}, m^i_{k,N},\lambda^i_{k,N})$ the solutions to the $N$--players game with viscosity $\nu=k\,\Id_d$, there hold
\begin{equation*}
\lim_{k\to0^+}\lim_{N\to+\infty}v^i_{k,N}=\lim_{N\to+\infty}\lim_{k\to0^+}v^i_{k,N} = v
\end{equation*}
in  $C^1_{loc}(\R^d)$ with second derivative converging uniformly in $\R^d$,
\begin{equation*}
~~\lim_{k\to0^+}\lim_{N\to+\infty}m^i_{k,N}=\lim_{N\to+\infty}\lim_{k\to0^+}m^i_{k,N} = m
\qquad\qquad \mbox{in distributional sense,}
\end{equation*}
\begin{equation*}
\!\!\!\!\!\!\!\!\!\!\!\!\!\!\!\!\!\!\!\!\!\!
\!\!\!\!\!\!\!\!\!\!\!\!\!\!\!\!\!
\lim_{k\to0^+}\lim_{N\to+\infty} \lambda^i_{k,N}=\lim_{N\to+\infty}\lim_{k\to0^+} \lambda^i_{k,N} = \lambda
\qquad\qquad ~~\mbox{ in $\R$,}
\end{equation*}
where 
%in the first case the limit is in $C^1_{loc}(\R^d)$ with second derivative converging uniformly in $\R^d$, in the second case the limit is in distributional sense, and 
$(v,m,\lambda)$ are given by
\begin{equation}\label{eq:vv_limit}
v(x)=\sqrt{r} \hat V(x-\hat \mu)+rAx\,,
\qquad\qquad
m=\delta_{\hat \mu}\,,
\qquad\qquad
\lambda = \hat F(0,\hat \mu)-\hat \mu^T\,{\hat V^2\over 2}\, \hat \mu\,,
\end{equation}
for $\hat V:=\sqrt{2\hat Q+r A^2}$, $\hat\mu:=(\hat V^2+\hat B)^{-1}\left(2\hat QH+\hat B\,\Delta\right)$ and $\hat F$ defined as in~\eqref{eq:v1_limit}.
\end{thm}

Finally, we study the limit when the cost for the control $r\to 0^+$, and thus large control can be chosen at cheap cost. Even if equations in~\eqref{eq:hjkfp} become singular when $r$ tends to zero, we can still use the formulas we found in the previous section to study the limit behavior.

\begin{thm}\label{thm:cc}  For $N\in\N$, consider $N$--players games of the form~\eqref{eq:sde}--\eqref{eq:cost_ergodic} such that {\bf (H1)}--{\bf (H3)} hold. Assume that~\eqref{eq:scale} holds as $N\to+\infty$ and that $\overline{\cal B}^\infty:=\hat Q+\,{\hat B\over 2}$ is invertible.

\n Then, the cheap control limit as $r\to 0^+$ and the mean field limit as $N\to+\infty$ commute.
Namely, denoting with $(v^i_{r,N}, m^i_{r,N},\lambda^i_{r,N})$ the solutions to the $N$--players game with control cost $R=r\,\Id_d$, there hold
\begin{equation*}
\lim_{r\to0^+}\lim_{N\to+\infty}v^i_{r,N}=\lim_{N\to+\infty}\lim_{r\to0^+}v^i_{r,N} = v
\end{equation*}
in  $C^1_{loc}(\R^d)$ with second derivative converging uniformly in $\R^d$,
\begin{equation*}
~~\lim_{r\to0^+}\lim_{N\to+\infty}m^i_{r,N}=\lim_{N\to+\infty}\lim_{r\to0^+}m^i_{r,N} = m
\qquad\qquad \mbox{in distributional sense,}
\end{equation*}
\begin{equation*}
\!\!\!\!\!\!\!\!\!\!\!\!\!\!\!\!\!\!\!\!\!\!
\!\!\!\!\!\!\!\!\!\!\!\!\!\!\!\!\!
\lim_{r\to0^+}\lim_{N\to+\infty} \lambda^i_{r,N}=\lim_{N\to+\infty}\lim_{r\to0^+} \lambda^i_{r,N} = \lambda
\qquad\qquad ~~\mbox{ in $\R$,}
\end{equation*}
where 
%in the first case the limit is in $C^1_{loc}(\R^d)$ with second derivative converging uniformly in $\R^d$, in the second case the limit is in distributional sense, and 
$(v,m,\lambda)$ are given by
\begin{equation}\label{eq:cc_limit}
v(x)\equiv 0\,,
\qquad\qquad
m=\delta_{\hat \mu}\,,
\qquad\qquad
\lambda = \hat F(0,\hat \mu)-\hat \mu^T\,{\hat V^2\over 2}\, \hat \mu\,,
\end{equation}
for $\hat V:=\sqrt{2\hat Q}$, $\hat\mu:=(\hat V^2+\hat B)^{-1}\left(2\hat QH+\hat B\,\Delta\right)$ and $\hat F$ defined as in~\eqref{eq:v1_limit}.
\end{thm}

\section{Technical proofs}\label{sec:tech_proof}

\n{\bf Proof of Theorem~\ref{thm:disc_Nplay}.} {\it Step 1.} By simply inserting the expressions of $v^i$ and $m^i$ into~\eqref{eq:hjkfp_disc}, one can transform the system of $2N$ equations into a system of equalities between quadratic forms to be satisfied for all $x\in\R^d$. Thus, by equating the coefficients of these quadratic forms, ~\eqref{eq:hjkfp_disc} reduces to algebraic relations~\eqref{eq:KFP_matrix_discount}--\eqref{eq:discount1} among the coefficients of $v^i$ and $m^i$.

\n It is now clear that if we show that there exists a unique solution in $\Symp$ to ARE in~\eqref{eq:discount1} for small $\ell$, then the existence and uniqueness part of the theorem would be proved. Indeed, the invertibility of ${\cal B}_\ell$ is equivalent to the existence and uniqueness of solutions for the linear system~\eqref{eq:discount1}, and once $\Sigma_\ell$ and $\mu_\ell$ are uniquely determined, conditions~\eqref{eq:KFP_matrix_discount} and~\eqref{eq:v0_disc} also give unique choices for $\Lambda_\ell,\rho_\ell, c^i_\ell$.

\n We therefore focus our attention on the ARE in~\eqref{eq:discount1}. By Proposition~\ref{prop:ARE}, solutions to~\eqref{eq:discount1} can be found as the $d$--dimensional invariant graph subspaces of the $2d\times 2d$ matrix ${\cal H}_\ell$ introduced in~\eqref{eq:discount_are_matrix}. Noticing that we have ${\cal A}_\ell\to 0$ and
${\cal Q}_\ell\to {\cal Q}=Q+r\,{A^2\over 2}$, as $\ell\to 0^+$, it is immediate to see that
$$
{\cal H}_\ell
\qquad\longrightarrow\qquad
{\cal H}:=\left(
\begin{array}{cc}
{\bf 0} & {\cal R}\\
{\cal Q} & {\bf 0}
\end{array}
\right)\,
$$
and that ${\cal H}$ has $d$ strictly positive and $d$ strictly negative eigenvalues. This latter property follows from the fact that $\lambda\in\mathrm{spec}({\cal H})$ if and only if $\lambda^2\in\mathrm{spec}({\cal RQ})$ and that Proposition~\ref{prop:matrices} implies $\mathrm{spec}({\cal RQ})\subset(0,+\infty)$ because both ${\cal R}$ and ${\cal Q}$ are positive definite. Therefore, all (possibly complex) eigenvalues of ${\cal H}_\ell$ will converge to some eigenvalue of ${\cal H}$, and there exists $\bar\ell>0$ small enough so that ${\cal H}_\ell$ has no non--zero purely imaginary eigenvalues when $\ell<\bar\ell$.
% purely imaginary eigenvalues of H_\ell have to converge to 0, but 0 is not eigenvalue of H so...
Propositions~\ref{prop:ARE}{\it (ii)} allows to conclude that ARE~\eqref{eq:discount1} admits symmetric solutions for $\ell<\bar\ell$.
%
% Note: solutions $\Sigma^i_{\ell^i}$ are real and symmetric, not only hermitian, because ${\cal H}_{\ell^i}$ is a real matrix!
%

\n Owing to Propositions~\ref{prop:ARE}{\it (iii)}, we also deduce that~\eqref{eq:discount1} has a unique symmetric solution $Y_\ell$ such that the eigenvalues of ${\cal A}_\ell+{\cal R}Y_\ell$ with non--zero real part are exactly the eigenvalues of ${\cal H}_\ell$ with positive real part. But
$$
{\cal A}_\ell+{\cal R}Y_\ell={k^2 r\over 2}\,\left(\,{\ell\over 2 k}\,\Id_d+Y_\ell\right)\,,
$$
is symmetric, so its eigenvalues are real and so are the ones of ${\cal H}_\ell$. Using~\eqref{eq:spec_ARE}, we obtain
$$
\mathrm{spec}({\cal A}_\ell+{\cal R}Y_\ell) = \mathrm{spec}({\cal H}_\ell)\cap (0,+\infty)\,.
$$
By setting $\delta:=\min\,\big\{\mathrm{spec}({\cal H})\cap (0,+\infty)\big\}>0$ and possibly reducing $\bar\ell$, we have for $\ell<\bar\ell$
$$
\min\,\big\{\mathrm{spec}({\cal H}_\ell)\cap (0,+\infty)\big\}~>\,{\delta\over 2}\,,
\qquad\qquad
\ell \,k\,r<\delta\,.
$$
Hence,
\begin{align*}
\min\,\mathrm{spec}(Y_\ell)~&=~{2\over k^2 r}\,\min\,\mathrm{spec}({\cal A}_\ell+{\cal R}Y_\ell)~-\,{\ell\over 2 k}\,
%&=~{2\over k^2 r}\left(\min\,\big\{\mathrm{spec}({\cal H}_\ell)\cap (0,+\infty)\big\}~-\,{\ell\, k\,r\over 4 }\right)\\
>\,{\delta\over 2k^2 r}\,>0\,,
\end{align*}
which implies $Y_\ell\in\Symp$ for $\ell<\bar\ell$. We claim that such a solution is also unique. Indeed, if any solution $Z_\ell\in\Symp$ exists with $Z_\ell\neq Y_\ell$, then $\mathrm{spec}({\cal A}_\ell+{\cal R}Z_\ell)=\,{k^2 r\over 2}\,\big(\mathrm{spec}(Z_\ell)+\,{\ell\over2 k}\big)\subseteq (0,+\infty)$ and this contradicts the characterization of $Y_\ell$ via~\eqref{eq:spec_ARE}. %This completes the first part of the proof.

\smallskip

\n{\it Step 2.} It remains to verify that affine feedback strategies $\overline{\alpha}^i(x)=(R^i)^{-1}\,{\nabla v^i(x)\over r}$ give a Nash equilibrium for the game~\eqref{eq:sde2}--\eqref{eq:cost_disc2}. Indeed, by applying Dynkin's formula,

\begin{align*}
\vatt&\Big[e^{-\ell T}v^i(X^i_T)-v^i(X^i_0)\Big]=
\vatt \Bigg[\int_0^T\!\!e^{-\ell s}\Big(\!-\ell v^i+\Tr(\nu\D^2 v^i)+(\nabla v^i)^T\!A\,x-(\nabla v^i)^T\alpha^i_s\Big)(X^i_s)ds \Bigg]\\
%&=\vatt \Bigg[\int_0^Te^{-\ell s}\Big(-\ell v^i+\Tr(\nu\D^2 v^i)+(\nabla v^i)^TA\,x-\big(R^{-1}\nabla v^i\big)^TR\,\alpha^i_s\Big)(X^i_s)\,ds \Bigg]\\
&\geq\vatt\Bigg[\int_0^T\!\!e^{-\ell s}\!\left(\!\!\Big(\!-\ell v^i+\Tr(\nu\D^2 v^i)+(\nabla v^i)^T\!A\,x- {(\nabla v^i)^T R^{-1}\nabla v^i\over 2}\,\Big)(X^i_s)-  \,{(\alpha^i_s)^TR\alpha^i_s\over 2}\,\right)\!ds\Bigg]\\
%&=\vatt \Bigg[\int_0^Te^{-\ell s}\Big(-\ell v^i(X^i_s)+\Tr(\nu\D^2 v^i(X^i_s))-H^i(X^i_s,\nabla v^i(X^i_s))\\
%&~~~~~~~~~~~~~~~~~~~~~~~~~~~~~~~~-(\alpha^i_s)^T \,{R\over 2}\,\alpha^i_s\Big)\,ds \Bigg]\\
&=-\vatt \Bigg[\int_0^T\!\!e^{-\ell s}\Big(f^i (X^i_s)+(\alpha^i_s)^T \,{R\over 2}\,\alpha^i_s\Big)ds \Bigg]
\end{align*}
with equality holding if $\alpha^i=\overline{\alpha}^i$. Since $\vatt\big[e^{-\ell T}v^i(X^i_T)\big]\to0$ as $T\to+\infty$, because the value function is quadratic and the strategies are admissible, we get 
\begin{align}\label{eq:fin_est}
v^i(X^i_0)&\leq \lim_{T\to+\infty}\vatt\left[\int_0^Te^{-\ell s}\Big(f^i (X^i_s)+(\alpha^i_s)^T \,{R\over 2}\,\alpha^i_s\Big)\,ds\right]\nonumber\\
&=~\vatt\left[\int_0^\infty e^{-\ell s}\Big(f^i (X^i_s)+(\alpha^i_s)^T \,{R\over 2}\,\alpha^i_s\Big)\,ds\right]
\end{align}
where we have used Lebesgue dominated convergence theorem in the last equality. 
% questa e' praticamente una successione di caratteristiche $\chi_{[0,T]}\to\chi_{[0,\infty[}$ che sono limitate dalla funzione
% |f(X_s)|+|R|/2 |\alpha_s|^2 che e' integrabile perche' funzione quadratica di X_s che e' L^2!
Noticing that equality holds only for $\alpha^i=\overline{\alpha}^i$, we can conclude that the cost corresponding to any unilateral change of strategy $\alpha^i$ (the r.h.s. of~\eqref{eq:fin_est}) is larger than the cost corresponding to $\overline{\alpha}^i$, and we have proved that $(\overline{\alpha}^1,\ldots,\overline{\alpha}^N)$ is a Nash equilibrium strategy.~~$\diamond$

\medskip

\n{\bf Proof of Theorem~\ref{thm:disc_LIMplay}}. 
{\it Step 1.} Proceeding as in the proof of Theorem~\ref{thm:disc_Nplay}, from imposing the expressions~\eqref{eq:disc_ansatz_limit} in~\eqref{eq:mfpde_disc}, we find that the coefficients $\Lambda,\Sigma,\rho,\mu$ satisfy the conditions~\eqref{eq:KFP_matrix_discount}--\eqref{eq:discount1} with $Q$ replaced by $\hat Q$ and with ${\cal B}_\ell$ replaced by ${\cal B}^\infty_\ell$. We can therefore repeat the arguments of the previous proof to show part {\it (a)}. In particular, we can assume that $\hat \ell>0$ is small enough to ensure that for $\ell<\hat\ell$ there hold $\mathrm{spec}({\cal H}_\ell)\subset\R$ and, setting $\ve^2:= \min\,\mathrm{spec}\left(\hat Q + r\,{A^2\over 2}\,\right) >0$,
\begin{equation}\label{eq:helper_disc}
\min\,\mathrm{spec}\left(\hat Q + r\,{A^2\over 2}\,-\,{\ell r\over2}\,A\right)>\,{\ve^2\over 2}\,,
\qquad\qquad
\ell\, k\, r< 2\ve\,.
\end{equation}

\smallskip

\n{\it Step 2.} Proceeding as in 
Theorem~4 in~\cite{BardiPriuli}, 
%Theorem~\textred{XXX} in~\cite{BardiPriuli}, 
it is easy to prove that $\hat B\geq 0$ is equivalent to the monotonicity of the operator $\hat V[m]$. Hence, we can repeat the arguments from~\cite{LL1,LL3} to show the uniqueness property {\it (b)}.

\smallskip

\n{\it Step 3.} As a preliminary step towards {\it (c)}, observe that $Q^N\to\hat Q$ as $N\to+\infty$ implies
$$
\min\,\mathrm{spec}\left(Q^N + r\,{A^2\over 2}\,-\,{\ell r\over2}\,A\right)>\,{1\over 2}\,\min\,\mathrm{spec}\left(\hat Q + r\,{A^2\over 2}\,-\,{\ell r\over2}\,A\right)> \,{\ve^2\over 4}\,>0\,,
$$
for $N$ large enough, where $\ve$ is the value introduced in step~1 and we have used~\eqref{eq:helper_disc} thanks to $\ell<\hat\ell$.
With the notations $Q^N_\ell:= Q^N + r\,{A^2\over 2}\,-\,{\ell r\over2}\,A$ and ${\cal H}^N_\ell:=\left(
\begin{array}{cc}
{\ell r k\over 4}\,\Id_d & \,{r k^2\over 2}\,\Id_d \\
Q^N_\ell & -\,{\ell r k\over 4}\,\Id_d
\end{array}
\right)
$,
we conclude that the matrix $Y^N_\ell\in\Sym$, solving the ARE corresponding to ${\cal H}^N_\ell$, satisfies
\begin{align*}
\min\,\mathrm{spec}(Y^N_\ell)&=\,{2\over k^2 r}~\min\,\{\mathrm{spec}({\cal H}^N_\ell)\cap(0,+\infty)\}-\,{\ell\over 2k}\,\\
&=\,{2\over k^2 r}~\sqrt{\min\,\mathrm{spec} \left({r k^2\over 2}\, Q^N_\ell\right)+\,{\ell^2 k^2 r^2\over 16}}-\,{\ell\over 2k}>\sqrt{2\over k^2 r}~\left(\,{\ve\over 2}\,-\,{\ell k r\over 4}\right)>0\,,
\end{align*}
because of the explicit expression of ${\cal H}^N_\ell$ in the second equality, and~\eqref{eq:helper_disc}. In particular, $Y^N_\ell$ is positive definite.
Observing that invertibility of ${\cal B}_\ell^\infty$also implies the invertibility of ${\cal B}_\ell$ for $N$ large enough, we conclude that~\eqref{eq:hjkfp_disc} admits a unique QG solution $(v^i_{\ell,N},m^i_{\ell,N})$ for large $N$.

\smallskip

\n{\it Step 4.} To pass to the limit as $N\to+\infty$, and complete the proof of part {\it (c)}, let us concentrate first on the sequence of the AREs in~\eqref{eq:discount1} as $N$ varies in $\N$. We can observe that ${\cal H}^N_\ell\to {\cal H}_\ell$ in~\eqref{eq:discount_are_matrix}, and that eigenvalues of ${\cal H}_\ell$ are real by our choice of $\ell$. Thus, the sequence of matrices $\Sigma^N_\ell$ solving~\eqref{eq:discount1} is bounded w.r.t. the norm $\|\cdot\|$ of the largest eigenvalue, defined in~\eqref{eq:max_spec_norm}, because
$$
\left\|\Sigma^N_\ell\right\|\,\leq\,{2\over k^2 r}\,\max\,\big\{\mathrm{spec}({\cal H}^N_\ell)\cap (0,+\infty)\big\}\, + \,{\ell\over 2k}
\leq\,{2\over k^2 r}\,\max\,\big\{\mathrm{spec}({\cal H})\cap (0,+\infty)\big\} +1 + \,{\ell\over 2k}\,,
$$
when $N$ is large enough.
There follows that $\Sigma^N_\ell$ has a converging subsequence $\Sigma^{N_m}_\ell$ whose limit $\overline{\Sigma}_\ell\in\Sym$ solves
$$
\,{r k^2\over 2}\,X^2+{\ell r k\over 2}\,X-\hat Q-r\,{A^2\over 2}\,+\ell\,r\,{A\over 2}\,=0\,,
$$
which is analogous to~\eqref{eq:discount1}, except for having $\hat Q$ in place of $Q$. If we could prove that $\overline{\Sigma}_\ell\in\Symp$, then we would have, by uniqueness in $\Symp$ of this limit ARE (which follows from {\bf (H1)} and Proposition~\ref{prop:ARE}{\it (iii)}), that $\overline{\Sigma}_\ell$ coincides with the matrix $\Sigma$ found in part {\it (a)} for the measure in~\eqref{eq:disc_ansatz_limit}. This additional property on $\overline{\Sigma}_\ell$ follows again by the continuity of the eigenvalues: we have seen in step~3 that for $N_m$ large enough we had
$$
\min\,\mathrm{spec}\left(\Sigma^{N_m}_\ell\right)\,>\sqrt{2\over k^2 r}~\left(\,{\ve\over 2}\,-\,{\ell k r\over 4}\right)>0\,,
$$
and this implies, as $N_m\to+\infty$, $\min\,\mathrm{spec}\big(\overline{\Sigma}_\ell\big)>0$, so that $\overline{\Sigma}_\ell\in\Symp$ and $\overline{\Sigma}_\ell=\Sigma$.

\n Now, we can pass to the limit $N\to +\infty$ also in the equation~\eqref{eq:discount1} for the average vector $\mu^N_\ell$: since ${\cal B}_\ell^\infty$ is invertible, we must have ${\cal B}_\ell$ invertible as well for $N$ large enough, so that
\begin{equation}
\mu^N_\ell={\cal B}_\ell^{-1}\Big(QH+(N-1)\,{B\over 2}\,\Delta\Big)
\qquad\longrightarrow\qquad
\mu=({\cal B}_\ell^\infty)^{-1}\Big(\hat QH+\,{\hat B\over 2}\,\Delta\Big)\,,
\end{equation}
i.e., $\mu^N_\ell$ converges to the average vector $\mu$ found in part {\it (a)}.
We conclude by observing that the previous convergence results for $\Sigma^N_\ell$ and $\mu^N_\ell$ allow to pass to the limit in~\eqref{eq:KFP_matrix_discount} and~\eqref{eq:v0_disc} as well, so to obtain the convergence of the value function.~~$\diamond$

\medskip

\n{\bf Proof of Theorem~\ref{thm:vd}.} {\it Step 1.} We start by proving that, when passing to the limit as $\ell\to 0^+$ in the discounted $N$--person game, the QG solution given in Theorem~\ref{thm:disc_Nplay} converges to the QG solution of the $N$--person game~\eqref{eq:sde}--\eqref{eq:cost_ergodic} given in section~\ref{sec:LQG}.

\n Let $\bar N\in\N$ be fixed large enough so that the matrix ${\cal B}^N= Q+\,{A^TRA\over 2}\,+(N-1)~{B\over 2}$ is invertible for $N\geq \bar N$ (compared to section~\ref{sec:LQG}, we added a superscript $N$ in the notation to stress its dependence on the number of players). For any $N\geq \bar N$, let us consider the discounted $N$--players game~\eqref{eq:sde2}--\eqref{eq:cost_disc2} and let $\bar\ell>0$ be the value found in Theorem~\ref{thm:disc_Nplay}. Since ${\cal B}^N_\ell:= Q+r\,{A^2\over 2}\,-\ell\,r\,{A\over 2}\,+(N-1)~{B\over 2}$ converges to ${\cal B}^N$ as $\ell\to 0^+$, it is not restrictive to assume that $\bar\ell$ is small enough to have ${\cal B}^N_\ell$ invertible for $\ell<\bar\ell$.

\n First, we focus our attention on the ARE in~\eqref{eq:discount1} and we fix any sequence $\ell_n\to 0^+$ with $\ell_n<\bar\ell$. By proceeding as in the proof of Theorem~\ref{thm:disc_LIMplay} above, we obtain that the sequence $\Sigma_{\ell_n}$ of solutions of~\eqref{eq:discount1} in $\Symp$ is bounded, and that any convergent subsequence has limit belonging to $\Symp$ and solving the ARE in~\eqref{eq:KFP_matrix}. Therefore, by uniqueness, we conclude that $\Sigma_{\ell_n}$ converges to the solution $\Sigma$ of~\eqref{eq:KFP_matrix} found in Theorem~2 of~\cite{BardiPriuli}.

\n By passing to the limit $\ell\to 0^+$ also in the equation for the average vector $\mu_\ell$ in~\eqref{eq:discount1}, we obtain
\begin{equation}
\mu_\ell=({\cal B}_\ell^N)^{-1}\Big(QH+(N-1)\,{B\over 2}\,\Delta\Big)
\qquad\longrightarrow\qquad
\mu=({\cal B}^N)^{-1}\Big(QH+(N-1)\,{B\over 2}\,\Delta\Big)\,,
\end{equation}
i.e., $\mu_\ell$ converges to the average vector $\mu$ found in~\eqref{eq:KFP_matrix}. In turn, $\Sigma_\ell\to\Sigma$ and $\mu_\ell\to\mu$ together with~\eqref{eq:KFP_matrix_discount}, imply $\Lambda_\ell\to\Lambda$ and $\rho_\ell\to\rho$ to the coefficients $\Lambda,\rho$ in~\eqref{eq:KFP_matrix}. 

\n Finally, from~\eqref{eq:v0_disc} we deduce easily that $c^i\to+\infty$, but also $\ell\, c^i\to\lambda^i$ and $\ell v^i_\ell(x)\to\lambda^i$ with $\lambda^i=F^i(\Sigma,\mu)+\Tr(\nu R\nu\Sigma+\nu R A)-\mu^T\,{\Sigma \nu R\nu \Sigma \over 2}\, \mu$ and $F^i$ given by~\eqref{eq:v1}, as in section~\ref{sec:LQG}.
Thus, we conclude
$$
v^i_\ell(x)-v^i_\ell(0)=x^T\,{\Lambda_\ell\over 2}\,x+\rho_\ell x~~ \longrightarrow ~~
x^T\,{\Lambda\over 2}\,x+\rho x=v^i(x)\,,
$$
recovering the expected value function of the problem with ergodic cost.

\smallskip

\n{\it Step 2.} Now we study the limit as $\ell\to 0^+$ of the mean field system~\eqref{eq:mfpde_disc}, and we fix $\tilde\ell>0$ small enough to have that $\ell<\tilde\ell$ implies invertibility of the matrix ${\cal B}_\ell^\infty$, defined in Theorem~\ref{thm:disc_LIMplay}.

\n For games with $\ell<\tilde\ell$, the part of step~1 about solutions of the ARE can be repeated, provided we replace $Q$ with $\hat Q$ in the various formulas derived from~\eqref{eq:mfpde_disc}. Namely, we can prove that the positive definite solutions $\hat\Sigma_\ell$ converge, as $\ell\to 0^+$, to the matrix $\Sigma\in\Symp$ which solves ARE in~\eqref{eq:KFP_matrix_limit}. Then, by passing to the limit in the equation for the average $\mu_\ell$ in~\eqref{eq:KFP_matrix_limit}, we obtain
\begin{equation}
\mu_\ell=({\cal B}_\ell^\infty)^{-1}\Big(\hat QH+\,{\hat B\over 2}\,\Delta\Big)
\qquad\longrightarrow\qquad
\mu=({\cal B}^\infty)^{-1}\Big(\hat QH+\,{\hat B\over 2}\,\Delta\Big)\,.
\end{equation}
The remaining coefficients converge like in step~1. In particular, $\ell v_\ell\to\lambda=\hat F(\Sigma,\mu)+\Tr(\nu R\nu\Sigma+\nu R A)-\mu^T\,{\Sigma \nu R\nu \Sigma \over 2}\, \mu$, with $\hat F$ given by~\eqref{eq:v1_limit}, as in section~\ref{sec:LQG}.

\smallskip

\n{\it Step 3.} By combining step~1 with the result on the mean field system~\eqref{eq:mfpde} in section~\ref{sec:LQG}, we obtain that
$$
\lim_{N\to+\infty}\lim_{\ell\to0^+}\Big[v^i_{\ell,N}-v^i_{\ell,N}(0)\Big]=v\,,
\qquad
\lim_{N\to+\infty}\lim_{\ell\to0^+}m^i_{\ell,N}=m\,,
\qquad
\lim_{N\to+\infty}\lim_{\ell\to0^+}\ell v^i_{\ell,N}=\lambda\,,
$$
in the appropriate topologies. 
Now, if we denote with $\hat\ell>0$ the minimum between the value found in Theorem~\ref{thm:disc_LIMplay}{\it (a)} and the value $\tilde \ell$ in step~2, for $\ell<\hat\ell$ the $N$--players game~\eqref{eq:sde2}--\eqref{eq:cost_disc2} has QG solutions, for all $N\in\N$. Moreover, taking $N$ large enough so that ${\cal B}^N_\ell$ is invertible (because it converges to the invertible matrix ${\cal B}_\ell^\infty$, as $N\to+\infty$), the QG solution is unique and it converges as $N\to+\infty$ to the solution of~\eqref{eq:mfpde_disc}, by Theorem~\ref{thm:disc_LIMplay}{\it (c)}. Thus, owing to step~2,
$$
\lim_{\ell\to0^+}\lim_{N\to+\infty}\Big[v^i_{\ell,N}-v^i_{\ell,N}(0)\Big]=v\,,
\qquad
\lim_{\ell\to0^+}\lim_{N\to+\infty}m^i_{\ell,N}=m\,,
\qquad
\lim_{\ell\to0^+}\lim_{N\to+\infty}\ell v^i_{\ell,N}=\lambda\,,
$$
so that~\eqref{eq:commute_disc1}--\eqref{eq:commute_disc3} hold and this concludes the proof.~~$\diamond$

\medskip

\n {\bf Proof of Theorem~\ref{thm:vv}.} {\it Step 1.} Using again the notation ${\cal B}^N= Q+\,{A^TRA\over 2}\,+(N-1)~{B\over 2}$, the convergence ${\cal B}^N\to{\cal B}^\infty$ as $N\to+\infty$ implies that there exists $\bar N\in\N$ such that ${\cal B}^N$ is invertible for $N\geq\bar N$. We then fix $k>0$, $N\geq \bar N$ and consider an $N$--players game satisfying assumptions {\bf (H1)}--{\bf (H3)} with $\nu=k\,\Id_d$. Instead of QG solutions of the form~\eqref{eq:ansatz_nearly_id}, we look for solutions to the HJB--KFP system~\eqref{eq:hjkfp} satisfying
\begin{equation}\label{eq:ansatz_h3}
v^i(x)=x^T\,{\Lambda\over 2}\, x+\rho x\,,\qquad\qquad
m^i(x)=\gamma\exp\left\{-\,{1\over 2} (x-\mu)^T\,{V\over k\,\sqrt{r}}(x-\mu)\right\}\,,
\end{equation}
for suitable matrices $V\in\Symp$, $\Lambda\in\Mat_{d\times d}(\R)$ and vectors $\mu,\rho\in\R^d$, which are the same for all the players. Here, $\gamma$ is a normalization constant explicitly given by $(2\pi)^{-d/2}\sqrt{\mathrm{det}(V^{-1})}$. 

\n By plugging these expressions into system~\eqref{eq:hjkfp}, or by setting $V= k\,\sqrt{r}\,\Sigma$ in the proof of Theorem~2 of~\cite{BardiPriuli}, one finds that the coefficients $\Lambda,V,\rho,\mu$ must satisfy
\begin{equation}\label{eq:diff_var1}
V^2=2Q^N+r\,A^2\,,
\quad
(V^2+(N-1)B^N)\,\mu=P\,,
\quad
\Lambda=\sqrt{r}\big(V+\sqrt{r}\,A\big)\,,
\quad \rho=-\sqrt{r}\,V\mu\,,
\end{equation}
%with $V^2+(N-1)B$ being invertible thanks to the assumption $B\geq 0$, 
where $P:=\left(2Q^NH+(N-1)B^N\Delta\right)$, and
\begin{equation}\label{eq:diff_var4}
\lambda^i=F^i\left({V\over k\sqrt{r}},\mu\right)-\mu^T\,{V^2\over 2}\, \mu+k \sqrt{r}\,\Tr(V+\sqrt{r}A)\,,
\end{equation}
with $F^i$ as in~\eqref{eq:v1}.
%and that all matrices $\Lambda^i$ are symmetric thanks to {\bf (H3)}.\\
It is immediate to check that, under our assumptions, the first two equations in~\eqref{eq:diff_var1} admit unique solution in $\Symp$ and $\R^d$, respectively, given by
\begin{equation}\label{eq:Vmu_vanish_visc}
V=\sqrt{2Q^N+r A^2}\,,
\qquad\qquad
\mu=(V^2+(N-1)B^N)^{-1}\left(2Q^NH+(N-1)B^N\Delta\right)\,,
\end{equation}
%which are well defined because $2Q+r A^2\in\Symp$ by {\bf (H2)} and because $B\geq 0$. 
Since~\eqref{eq:diff_var1} do not depend on $k$, the same is true for the value function $v^i$ and for the mean vector $\mu$. Only the value $\lambda^i$ in~\eqref{eq:diff_var4} is modified by a change of $k$.
Passing finally to the limit as $k\to 0^+$ in~\eqref{eq:diff_var1}--\eqref{eq:diff_var4}, we conclude
$$
v^i_{k,N}(x)\to\sqrt{r} V(x-\mu)+rAx\,,
\quad
m^i_{k,N}={\cal N}\left(\mu,{V_k\over k\sqrt{r}}\right)\to \delta_\mu\,,
\quad
\lambda^i_{k,N}\to F^i(0,\mu)-\mu^T\,{V^2\over 2}\, \mu\,,
$$
in the correct topologies,  with $V$ and $\mu$ given by~\eqref{eq:Vmu_vanish_visc}.

\smallskip

\n{\it Step 2.} Analogous computations can be performed for the mean field equations~\eqref{eq:mfpde}. In this case, the expression~\eqref{eq:vv_limit} for the value function remains valid as $k\to 0^+$, and it is easy to verify that $m={\cal N}\left(\hat \mu, {\hat V\over k\sqrt{r}}\right)\to \delta_{\hat\mu}$ in distributional sense. Since we also have
$$
\lambda=\hat F\left({\hat V\over k\sqrt{r}},\hat \mu\right)-\hat \mu^T\,{\hat V^2\over 2}\, \hat \mu+k \sqrt{r}\,\Tr(\hat V+\sqrt{r}A)
~\longrightarrow~
\hat F(0,\hat \mu)-\hat \mu^T\,{\hat V^2\over 2}\, \hat \mu\,,
$$
it is enough to pass to the limit as $N\to\infty$ in the formulas~\eqref{eq:Vmu_vanish_visc} to complete the proof.
%%
%%\medskip
%%
%%\n{\bf 3.} Notice that, thanks to the explicit formulas in~\eqref{eq:Vmu_vanish_visc}, $V\to\hat V$ and $\mu\to\hat\mu$ as $N\to+\infty$. Hence, the expressions found in {\bf 1.} satisfy
%%$$
%%\sqrt{r} V(x-\mu)+rAx ~~\to ~~\sqrt{r} \hat V(x-\hat \mu)+rAx\,,
%%\qquad\qquad\qquad
%%\delta_\mu ~~\to ~~\delta_{\hat \mu}\,,
%%$$$$
%%F^i(0,\mu)-\mu^T\,{V^2\over 2}\, \mu ~~\to ~~ \hat F(0,\hat\mu)-\hat\mu^T\,{\hat V^2\over 2}\, \hat\mu\,,
%%$$
%%as $N\to+\infty$, and this completes the proof.
~~$\diamond$

\medskip

\begin{rem}\label{rem:extenstion} Since~\eqref{eq:diff_var1} is independent from the matrix $\nu$, the assumption  {\bf (H3)} in Theorem~\ref{thm:vv} can be replaced by the following, slightly more general, one:
\begin{description}
\item{{\bf (H3$'$)}} The matrix $A$ is symmetric, there exist a constant $r>0$ such that $R=r\,\Id_d$, and there holds $\nu = k\,\bar\nu$ for a constant $k>0$ and a matrix $\bar\nu$ such that $\bar\nu \sqrt{2Q+r A^2}\in\Sym$.
\end{description}
\end{rem}

\medskip

\n {\bf Proof of Theorem~\ref{thm:cc}.} {\it Step 1.} The convergence $\overline{\cal B}^N:=Q^N+(N-1)\,{B^N\over 2}\to\overline{\cal B}^\infty$ as $N\to+\infty$ implies that there exists $\bar N\in\N$ such that $\overline{\cal B}^N$ is invertible for $N\geq\bar N$. Fixed $N\geq \bar N$, we thus consider the $N$--players game~\eqref{eq:sde}--\eqref{eq:cost_ergodic} and let $r>0$ small enough so that ${\cal B}^N:=Q^N+r\,{A^2\over 2}+(N-1)\,{B^N\over 2}$ is invertible too.

\n By looking for solutions of the form~\eqref{eq:ansatz_h3} in the HJB--KFP system with cost $R=r\,\Id_d$, we find that matrices $\Lambda,V$ and vectors $\rho,\mu$ satisfy again~\eqref{eq:diff_var1}--\eqref{eq:diff_var4}. 

\n When $r\to 0^+$, we claim that the value functions $v^i\to 0$ uniformly on compact sets. Indeed, for fixed $r>0$, the ARE in~\eqref{eq:diff_var1} admits a unique solution $V_r\in\Symp$, given by
$V_r:=\sqrt{2Q^N+r A^2}$.
As $r\to 0^+$, we thus have $V_r\to\overline{V}=\sqrt{2Q^N}$. %, which is the unique symmetric positive definite solution to the limit ARE $\overline{V}^2=2Q$. 
Similarly, by passing to the limit in the equation~\eqref{eq:diff_var1} for $\mu$, we obtain
$$
\mu_r:=(V_r^2+(N-1)B^N)^{-1}P
~~\longrightarrow~~
\overline{\mu}:=(\overline{V}^2+(N-1)B^N)^{-1}P\,,
$$
which in turn implies $\Lambda\to 0$ and $\rho\to 0$. It is also simple to verify that the measures $m^i$ converge in distributional sense to a Dirac delta $\delta_{\overline{\mu}}$, centered at $\overline{\mu}$, and that
$\lambda^i_{r,N}
\longrightarrow
F^i(0,\overline{\mu})-\overline{\mu}^T\,{\overline{V}^2\over 2}\, \overline{\mu}$,
as in the deterministic limit.

\smallskip

\n {\it Step 2.} Fixed $r>0$ small enough to have invertibility of the matrix ${\cal B}^\infty:=\hat Q+r\,{A^2\over 2}+\,{\hat B\over 2}$, we repeat the argument used in step~1 for the mean field equations~\eqref{eq:mfpde}. In this case, we get
$$
v_r(x) = \sqrt{r} \hat V_r(x-\mu)+rAx\,,
\qquad\qquad
m_r= {\cal N}\left(\hat\mu_r,{\hat V_r\over k\sqrt{r}}\right)\,,
$$
with $\hat V_r:=\sqrt{2\hat Q+r A^2}$ and $\hat\mu_r:=(\hat V_r^2+\hat B)^{-1}\left(2\hat QH+\hat B\,\Delta\right)$.
Then, as $r\to 0^+$ we obtain the convergence of $v_r$ and $m_r$ to the value function and the measure in~\eqref{eq:cc_limit}.
From 
$$
\lambda_r=\hat F\left({\hat V_r\over k\sqrt{r}},\hat \mu_r\right)-\hat \mu_r^T\,{\hat V_r^2\over 2}\, \hat \mu_r+k \sqrt{r}\,\Tr(\hat V_r+\sqrt{r}A)
~\longrightarrow~
\hat F(0,\hat \mu)-\hat \mu^T\,{\hat V^2\over 2}\, \hat \mu\,.
$$
also the convergence of $\lambda_r$ follows. Finally, by passing to the limit as $N\to+\infty$ in the formulas in step~1, it is immediate to prove $\overline{V}\to\hat V$ and $\overline{\mu}\to\hat\mu$, whence the conclusion follows. 
%%
%%\medskip
%%
%%\n{\bf 3.} Observe that, thanks to their explicit formulas in {\bf 1.}, it is immediate to prove $\overline{V}\to\hat V$ and $\overline{\mu}\to\hat\mu$ as $N\to+\infty$. Hence, as $N\to+\infty$ we find that
%%$$
%%\delta_{\overline{\mu}} ~~ \longrightarrow ~~\delta_{\hat \mu}\,,
%%\qquad\qquad\qquad
%%F^i(0,\overline{\mu})-\overline{\mu}^T\,{\overline{V}^2\over 2}\, \overline{\mu} ~~ \longrightarrow ~~ \hat F(0,\hat\mu)-\hat\mu^T\,{\hat V^2\over 2}\, \hat\mu\,,
%%$$
%%and this completes the proof.
~~$\diamond$

\section{Extensions and open problems}\label{sec:conclusion}

\subsection{Games with $N$ different players}

For sake of notational simplicity, in this work we have focused our attention on games satisfying {\bf (H2)}, i.e. games with nearly identical players (see also Definition~4.1 in~\cite{BardiPriuli}). 
%\textred{Definition~XXX} in~\cite{BardiPriuli}). 
However, some of the results we have presented admit a straightforward generalization to games whose cost for the player's state~\eqref{eq:quad_cost} consists of more general matrix coefficients $Q^i_{jk}$. The most interesting extension is probably the characterization of affine Nash equilibria strategies for general games with discounted cost~\eqref{eq:cost_disc}.
We replace assumptions {\bf (H1)}--{\bf (H3)} with the following
\begin{description}
\item{\bf (H)} The matrix $\sigma$ is invertible and the matrix $R$ belongs to $\Symp$. Moreover, for all $i\in\{1,\ldots,N\}$, let assume that matrices $Q^i$ are symmetric, that the block $Q^i_{ii}\in\Symp$ and that $(A,\nu,R,Q^i_{ii})$ satisfy the Riccati--Sylvester property in the sense of Definition~\ref{def:RSprop}.
\end{description}
In this case system~\eqref{eq:hjkfp_disc_intro} takes the same form as~\eqref{eq:hjkfp}, but with $\ell v^i$ in place of $\lambda^i$ in the first equation for each player. Then, we can prove the following result.

\begin{thm}\label{thm:disc_Nplay2} Under assumption {\bf (H)}, there exists $\bar\ell>0$ such that for $\ell<\bar\ell$ the system of HJB--KFP equations for the game admits a unique solution $(v^i,m^i)$ of the form
$$
v^i(x)=x^T\,{\Lambda^i\over 2}\, x+\rho^i x+c^i\,,\qquad\qquad
m^i(x)= {\cal N}(\mu^i,(\Sigma^i)^{-1})\,,
$$
for suitable symmetric matrices $\Lambda^i,\Sigma^i$, with $\Sigma^i$ positive definite, vectors $\mu^i,\rho^i$ and numbers $c^i$, if and only if the $Nd \times Nd$ matrix
$$
\widetilde{\cal B}:=\big(\widetilde{\cal B}_{\alpha\beta}\big)_{\alpha,\beta=1,\ldots,N}
%$$
\qquad\qquad
%$$
\widetilde{\cal B}_{\alpha\beta}:= Q^\alpha_{\alpha\beta}+\delta_{\alpha\beta}\left(\,{A^T R A\over 2}\,-\ell\,{R A\over 2}\right)\,\in\Mat_{d\times d}(\R)\,.
$$
is invertible, $\delta_{\alpha\beta}$ being the Kronecker delta.
Moreover, the affine feedbacks
$\overline{\alpha}^i(x)=R^{-1}\nabla v^i(x)$, for $x\in\R^d$ and $i=1,\ldots,N$,
provide a Nash equilibrium strategy for the game~\eqref{eq:sde}--\eqref{eq:cost_disc}, for all initial states $X\in\R^{Nd}$.
\end{thm}

\n The proof proceeds along the same lines of the one for Theorem~\ref{thm:disc_Nplay}, and it is therefore omitted. Further extensions to games having matrices $A^i,\sigma^i,R^i$ and discount factors $\ell^i$ also depending on the players just require changes in the corresponding notations.

\subsection{Games not satisfying {\bf (H3)}}

Comparing the results presented in this paper with the ones in~\cite{BardiPriuli}, one might easily wonder why the assumption {\bf (H3)} is here imposed on some matrix coefficients in the dynamics and the cost. The answer is related to the algebraic Riccati equations whose solutions give the (inverse of the) covariance matrix of the desired Gaussian measure. Indeed, in the case of deterministic and cheap control limits, a large part of the manipulations done on the system~\eqref{eq:hjkfp} of HJB--KFP equations can still be repeated for games not satisfying {\bf (H3)}. By searching for solutions of the form
\begin{equation}\label{eq:ansatz_gen}
v^i(x)=x^T\,{\Lambda\over 2}\, x+\rho x\,,\qquad\qquad
m^i(x)=\gamma\exp\left\{-\,{1\over 2} (x-\mu)^T\,\nu^{-1}R^{-1/2}V\,(x-\mu)\right\}\,,
\end{equation}
one finds relations similar to~\eqref{eq:diff_var1} and, in particular, we have that $V$ must solve
$V^TV=2Q+A^TRA$.
However, in this context we are not searching for solutions $V\in\Sym$ anymore, but for a $V$ which makes $\Sigma:=\nu^{-1}R^{-1/2}V\in\Symp$.

For any fixed choice of the matrices $\nu$ and $R$, the existence and uniqueness result in Theorem~2 of~\cite{BardiPriuli} (see also section~\ref{sec:LQG}) allows to prove that a unique $V$ with the required properties exists, and thus that a unique QG solution to~\eqref{eq:hjkfp} exists of the form~\eqref{eq:ansatz_gen}, at least when $(A,\nu,R,Q)$ satisfy the Riccati--Sylvester property in the sense of Definition~\ref{def:RSprop} and ${\cal B}= Q+\,{A^TRA\over 2}\,+(N-1)~{B\over 2}$ is invertible. 
The problems for these more general games arise when we try to pass to the limit: Indeed, except for the simple extension mentioned in Remark~\ref{rem:extenstion}, it is not clear whether the sequences of solutions converge, either as $\nu\to 0$ or as $R\to 0$, to a specific limit matrix $\overline{V}$ among the many solutions of the limit ARE equation which is, respectively,
$$
\overline{V}^T\overline{V}=2Q+A^TRA\,,
\qquad\qquad
\overline{V}^T\overline{V}=2Q\,.
$$
Analogous issues are found when studying the limits of mean field equations~\eqref{eq:mfpde}.

For the vanishing discount limit there is an additional difficulty, because it is not clear whether symmetric positive definite solutions to the ARE in~\eqref{eq:discount1} exist when {\bf (H3)} is not satisfied. Indeed, the matrix ${\cal A}_\ell+{\cal R}Y_\ell=\,{\nu R\over 2}\,\left({\ell\over 2}\,\Id_d+\nu Y_\ell\right)$, with $Y_\ell\in\Sym$ given by Proposition~\ref{prop:ARE}{\it (iii)}, might be not symmetric and have complex eigenvalues if we do not assume {\bf (H3)}. In this case, it does not seem possible to generally deduce $Y_\ell>0$ from the estimates on the real part of eigenvalues of ${\cal A}_\ell+{\cal R}Y_\ell$. 

In our opinion, new results on the algebraic Riccati equations would be necessary to extend our analysis to more general games, but such extensions are beyond the scope of this work.

%%In conclusion, we restricted our investigation to games satisfying {\bf (H3)} because most results on ARE do not deal with stability (resp. instability) of the solution itself. Indeed, they have been developed in connection with control problems for affine systems  $\dot x=Ax+Bu$, where the goal is the stability (resp. instability) of the closed loop system
%%$$
%%\dot x = (A+BF)x
%%$$
%%corresponding to a linear feedback $u(x)=Fx$, with $F$ solving  a suitable algebraic Riccati equation, and do not provide enough information for the case we are considering here.
%%In other words, some new results on the algebraic Riccati equations would be necessary to extend our analysis to more general games, but these extensions are beyond the scope of this work.
%%
\subsection{Comparison with previous works on infinite horizon games with discounted cost}\label{sec:compare_caines}

For $N$--person infinite horizon games with discounted costs there is a rich literature (see~\cite{HCM07ieee,LZ} and references therein). Typically, the games considered have a dynamics
\begin{equation}\label{eq:caines_sde}
dX^i_t=\left({\bf A}\,X^i_t+{\bf B}\,\alpha^i_t\right)dt+{\bf D}\,dW^i_t\,,
\qquad\qquad
i=1,\ldots,N\,,
\end{equation}
and cost
\begin{equation}\label{eq:caines_cost}
J^i(X,\alpha^1,\ldots,\alpha^N):= \vatt\left[\int_0^{+\infty}e^{-\ell t}\left(\,{(\alpha_t^i)^T\,{\bf R}\,\alpha_t^i\over 2}\,+(X^i_t-\Xi_t)^T {\bf Q}\,(X_t-\Xi_t)\right)\,dt\right]\,,
\end{equation}
where ${\bf A,B,D,R,Q}\in\Mat_{d\times d}(\R)$ suitable matrices, ${\bf R}>0$ and ${\bf Q}\geq 0$, and where
$$
\Xi_t := \Gamma \cdot\left(\,{1\over N}\,\sum_{j=1}^N X^j_t\right) +\eta\,\in\R^d \,,
\qquad\qquad \Gamma\in\Mat_{d\times d}(\R)\,,~\eta\in\R^d\,,
$$
with the term $1/N\sum_{j=1}^N X^j_t$ representing a sort of average position among the agents (referred to as the ``mean field term'' of the game). The typical result for these games is that the solution of a suitable ``mean field system'' of ODEs, obtained by formally passing to the limit as $N\to+\infty$ in the HJB equation for~\eqref{eq:caines_sde}--\eqref{eq:caines_cost} and replacing the mean field term in the cost with a suitable deterministic function, provides an approximate Nash equilibria for the game~\eqref{eq:caines_sde}--\eqref{eq:caines_cost}. Namely, the feedback strategy corresponding to such a solution is an $\ve$--Nash equilibrium strategy with $\ve={\cal O}\left({1\over N}\right)$.

\smallskip

Now observe that the second term in the cost~\eqref{eq:caines_cost} can be rewritten in the form~\eqref{eq:quad_cost}, by choosing $\overline{X}^i_i=\eta$, $\overline{X}^j_i=0$ for $j\neq i$, and
$$
Q^i_{ii}=\left(\Id_d-\,{\Gamma^T\over N}\right){\bf Q}\left(\Id_d-\,{\Gamma\over N}\right)\,,
%$$
\qquad
%$$
Q^i_{ij}=-\left(\Id_d-\,{\Gamma^T\over N}\right){\bf Q}\,{\Gamma\over N}\,,
%$$
\qquad
%$$
Q^i_{jk}=\,{\Gamma^T\over N}\,{\bf Q}\,{\Gamma\over N}\,,
$$
and that {\bf (H1)} is satisfied whenever ${\bf Q}>0$, since $\Id_d-\,{\Gamma^T\over N}$ is always invertible for $N$ large enough. Therefore, games~\eqref{eq:caines_sde}--\eqref{eq:caines_cost} are very similar to the ones we considered in Theorem~\ref{thm:disc_Nplay}. The main difference is that the cost $J^i$ in~\eqref{eq:caines_cost} depends on other players directly through their state $X^j$, while in~\eqref{eq:cost_disc} the dependence is present only through their asymptotic distribution $m^j$ in the environment. The novelty in our results is that, thanks to the particular form of the cost, we are able to characterize \emph{exact} Nash equilibria for the discounted game (at least for small values of the discount factor $\ell$) and not only of $\ve$--approximate ones. Moreover, we prove rigorously the convergence of such Nash equilibria to the solutions of the mean field game. The analogous study in the case of games with cost~\eqref{eq:caines_cost} is still an open problem, to our knowledge.

\end{document}